\ifx\shlhetal\undefinedcontrolsequence\let\shlhetal\relax\fi

\documentclass[11pt]{amsart}

\usepackage{amsmath}
\usepackage{amssymb}

\newtheorem{theorem}{Theorem}[section]
\newtheorem{claim}[theorem]{Claim}

\newtheorem{lemma}[theorem]{Lemma}

\newtheorem{proposition}[theorem]{Proposition}
\newtheorem{corollary}[theorem]{Corollary}

\theoremstyle{definition}
\newtheorem{definition}[theorem]{Definition}

\theoremstyle{remark}

\newtheorem{hypothesis}[theorem]{Hypothesis}

\newcount\skewfactor
\def\mathunderaccent#1#2 {\let\theaccent#1\skewfactor#2
\mathpalette\putaccentunder}
\def\putaccentunder#1#2{\oalign{$#1#2$\crcr\hidewidth
\vbox to.2ex{\hbox{$#1\skew\skewfactor\theaccent{}$}\vss}\hidewidth}}

\def\smallbox#1{\leavevmode\thinspace\hbox{\vrule\vtop{\vbox
   {\hrule\kern1pt\hbox{\vphantom{\tt/}\thinspace{\tt#1}\thinspace}}
   \kern1pt\hrule}\vrule}\thinspace}

\newcommand{\cf}{{\rm cf}}

\newcommand{\then}{{\underline{then}}}
\newcommand{\Then}{{\underline{Then}}}

\def\qedref#1{$\qed_{\reforiginal{#1}}$}

\title{pity on $\lambda$}
\author{Shimon Garti}
\address{Institute of Mathematics
 The Hebrew University of Jerusalem
 Jerusalem 91904, Israel}
\email{shimon.garty@mail.huji.ac.il}


\begin{document}
\let\labeloriginal\label
\let\reforiginal\ref

\begin{abstract}
Assume $\lambda = \lambda^{<\lambda}$. We prove a generalization of Rothberger's theorem, namely $\mathfrak{p}_\lambda = \lambda^+ \Rightarrow \mathfrak{t}_\lambda = \lambda^+$. We conclude from this theorem that $\cf(\mathfrak{p}_\lambda) \neq \lambda$ for every $\lambda = \lambda^{<\lambda}$. We also say something about the discrepancy between $\mathfrak{p}_\lambda$ and $\mathfrak{t}_\lambda$.
\end{abstract}

\maketitle

\newpage

\section{introduction}

The celebrated theorem of Rothberger (see \cite{MR0004281}) asserts that if $\mathfrak{p} = \aleph_1$ then $\mathfrak{t} = \aleph_1$. $\mathfrak{p}$ and $\mathfrak{t}$ are two cardinal invariants on the continuum, both are uncountable. The relation $\mathfrak{p} \leq \mathfrak{t}$ is trivial, and the question whether $\mathfrak{p} < \mathfrak{t}$ is possible is a longstanding open problem in set theory. On the consequences of $\mathfrak{p} < \mathfrak{t}$ (if consistent), see \cite{MR2518968}. Nevertheless, if $\mathfrak{p} = \aleph_1$ then $\mathfrak{p} = \mathfrak{t}$.

A wonderful exposition to the theorem of Rothberger appears in \cite{MR776622}. It is explained there what is the difficulty in trying to stretch this theorem to cardinals above $\aleph_1$. It is proved there that the basic argument of Rothberger cannot be applied to the case of $\mathfrak{p} \geq \aleph_2$.

But there is an interesting case that carries Rothberger's argument even if $\mathfrak{p} > \aleph_1$, the case of singular cardinals with countable cofinality. Let us try to explain the idea. Rothberger proved that if $\mathfrak{p} = \aleph_1$ then \emph{every} $\mathfrak{p}$-family of size $\aleph_1$ can be converted into a tower of size $\aleph_1$. We shall see that if $\lambda > \cf(\lambda) = \aleph_0$ then $\mathfrak{p}$-families of size $\lambda$ (here, a $\mathfrak{p}$-family is a family of sets with the strong finite intersection property, and no infinite pseudo intersection; such a family does not exemplify the cardinal invariant $\mathfrak{p}$ of course) can be converted into towers of size $\lambda$. But such a tower yields a tower of size $\cf(\lambda)$, which is absurd.

The conclusion is that no $\mathfrak{p}$-family of size $\lambda$ exists. How can we use it? Well, we shall see that if a singular cardinal with countable cofinality separates between $\mathfrak{p}$ and $\mathfrak{t}$ then one can find (under some reasonable assumption) a $\mathfrak{p}$-family of size $\lambda$. We can conclude that under this reasonable assumption $\mathfrak{t} \leq \mathfrak{p}^{+n}$ for some $n \in \omega$.

Another way to exploit this idea, applies to $\mathfrak{p}_\lambda, \mathfrak{t}_\lambda$ for $\lambda > \aleph_0$. Here we are interested in the following generalization of Rothberger's theorem. Let $\lambda$ be an infinite cardinal. We can define $\mathfrak{p}_\lambda$ and $\mathfrak{t}_\lambda$ in a similar way to the original definition of $\mathfrak{p}$ and $\mathfrak{t}$. We will show that if $\mathfrak{p}_\lambda = \lambda^+$ then $\mathfrak{t}_\lambda = \lambda^+$, whenever $\lambda = \lambda^{<\lambda}$.

As a consequence we can derive a conclusion on the cofinality of $\mathfrak{p}_\lambda$. A straightforward argument shows that $\mathfrak{t}_\lambda$ is always a regular cardinal. It is known also that $\mathfrak{p}$ is regular, but the proof is more sophisticated. Moreover, it is not clear if this proof applies to $\mathfrak{p}_\lambda$ for higher $\lambda$-s. From the generalization to Rothberger's theorem we will be able to prove that the cofinality of $\mathfrak{p}_\lambda$ is at least $\lambda^+$. Another consequence is related to the distance between $\mathfrak{p}_\lambda$ and $\mathfrak{t}_\lambda$.

We try to adhere to the standard notation. We denote infinite cardinals by $\theta, \kappa, \lambda, \mu, \tau$. We use $\alpha, \beta, \gamma, \delta, \varepsilon, \zeta, \eta, \xi$ for ordinals. By $A \subseteq^*_\lambda B$ we mean that $|A \setminus B| < \lambda$, but in most cases we shall write $A \subseteq^* B$ and omit $\lambda$ since it will be clear from the context. the notation concerns topological invariants on the continuum is due to \cite{MR776622}.

\newpage

\section{Background}

\begin{hypothesis}
\label{rregularity}
Throughout the paper, $\aleph_0\leq\lambda=\lambda^{<\lambda}$.
\end{hypothesis}

\begin{definition}
\label{bbb}
The bounding number, $\mathfrak{b}_\lambda$.
\begin{enumerate}
\item [$(\aleph)$] For $f,g \in {^\lambda \lambda}$ we define $f \leq^* g$ iff $|\{ \alpha < \lambda : f(\alpha) > g(\alpha) \}| < \lambda$
\item [$(\beth)$] $B \subseteq {^\lambda \lambda}$ is unbounded if there is no $h \in {^\lambda \lambda}$ so that $f \leq^* h$ for every $f \in B$
\item [$(\gimel)$] $\mathfrak{b}_\lambda = {\rm min} \{|B| : B \subseteq {^\lambda \lambda}, B$ is unbounded $\}$
\end{enumerate}
\end{definition}

If $\lambda = \omega$ then the definition of $\mathfrak{b}_\lambda$ coincides with the common definition of the bounding number $\mathfrak{b}$. In a similar way we define $\mathfrak{d}_\lambda$ as the dominating number in ${^\lambda \lambda}$. A set $B \subseteq {^\lambda \lambda}$ is dominating when for every $f \in {^\lambda \lambda}$ there exists $h \in B$ so that $f \leq^* h$. $\mathfrak{d}_\lambda$ is the minimal cardinality of a dominating set. Clearly, $\mathfrak{b}_\lambda \leq \mathfrak{d}_\lambda$. Trying to follow the methods of proof in \cite{MR776622}, we define another version of $\mathfrak{b}_\lambda$ and prove that both definitions are equivalent.

\begin{definition}
\label{bbb1}
The strictly increasing bounding number, $\mathfrak{b}^1_\lambda$. \newline 
$\mathfrak{b}^1_\lambda = {\rm min} \{|B| : B \subseteq {^\lambda \lambda}, B$ is unbounded, the members of $B$ are strictly increasing functions, and $B$ is well ordered by $<^* \}$
\end{definition}

\begin{lemma}
\label{bbbisbbb1}
$\mathfrak{b}_\lambda = \mathfrak{b}^1_\lambda$.
\end{lemma}

\par \noindent \emph{Proof}. \newline 
$\mathfrak{b}^1_\lambda \geq \mathfrak{b}_\lambda$, since every $B$ which attests $\mathfrak{b}^1_\lambda$ is also an evidence to $\mathfrak{b}_\lambda$. Let us show that $\mathfrak{b}^1_\lambda \leq \mathfrak{b}_\lambda$. Choose $B = \{f_\eta : \eta < \mathfrak{b}_\lambda \}$, an evidence to $\mathfrak{b}_\lambda$. Let $S$ be the collection of the strictly increasing functions in ${^\lambda \lambda}$. Notice that $S$ is a dominating subset of ${^\lambda \lambda}$.

For every $\alpha < \mathfrak{b}_\lambda$ choose $f'_\alpha \in S$ which is strictly above $\{f'_\beta : \beta < \alpha\} \cup \{f_\alpha\}$ (exists, since the size of this set is below $\mathfrak{b}_\lambda \leq \mathfrak{d}_\lambda$). Set $B' = \{f'_\eta : \eta < \mathfrak{b}_\lambda \}$. Clearly, $|B'| = |B|$ and $B'$ is an evidence to $\mathfrak{b}^1_\lambda$, so we are done.

\hfill \qedref{bbbisbbb1}

The following definition is also a straightforward generalization of the usual notions, applied to $\lambda$:

\begin{definition}
\label{ppp}
The pseudointersection number, $\mathfrak{p}_\lambda$. \newline 
Let $\mathcal{F}$ be a family of subsets of $\lambda$:
\begin{enumerate}
\item [$(\aleph)$] $A$ is a pseudointersection of $\mathcal{F}$ if $A \subseteq^* F$ for every $F \in \mathcal{F}$
\item [$(\beth)$] $\mathcal{F}$ has the strong intersection property if for every $\mathcal{F'} \subseteq \mathcal{F}$ so that $|\mathcal{F'}| < \lambda$ we have $|\bigcap \mathcal{F'}| = \lambda$
\item [$(\gimel)$] $\mathfrak{p}_\lambda = {\rm min} \{|\mathcal{F}| : \mathcal{F} \subseteq [\lambda]^\lambda, \mathcal{F}$ has the strong intersection property, and there is no pseudointersection of $\mathcal{F}$ of size $\lambda \}$
\end{enumerate} 
\end{definition}

The cardinal invariant $\mathfrak{p}$ (i.e., $\mathfrak{p}_\omega$ in our notation) is always a regular cardinal. It is shown in section $3$ below that $\cf(\mathfrak{p}_\lambda) > \lambda$ when $\lambda = \lambda^{<\lambda}$. The last basic definition that we need is the following:

\begin{definition}
\label{ttt}
The towering number, $\mathfrak{t}_\lambda$.
\begin{enumerate}
\item [$(\aleph)$] $\mathcal{T} \subseteq [\lambda]^\lambda$ is a tower in $\lambda$, if $\mathcal{T}$ is well ordered by $\supseteq^*$, $\mathcal{T}$ has the strong intersection property, and there is no pseudointersection for $\mathcal{T}$ of size $\lambda$
\item [$(\beth)$] $\mathfrak{t}_\lambda = {\rm min} \{|\mathcal{T}| : \mathcal{T}$ is a tower in $\lambda\}$
\end{enumerate}
\end{definition}

Notice that $\mathfrak{p}_\lambda \leq \mathfrak{t}_\lambda$ since every tower has the strong intersection property.
The cardinal invariant $\mathfrak{t}$ on the continuum is always regular. This is generalized easily to the cardinal $\mathfrak{t}_\lambda$ for every $\lambda$:

\begin{lemma}
\label{ttlambdareg}
The regularity of $\mathfrak{t}_\lambda$. \newline 
Let $\lambda$ be an infinite cardinal. \then\ $\mathfrak{t}_\lambda$ is a regular cardinal.
\end{lemma}

\par \noindent \emph{Proof}. \newline 
Denote $\mathfrak{t}_\lambda$ by $\mu$. Assume toward contradiction that $\mu > \cf(\mu) = \kappa$. Choose an evidence to the fact that $\mathfrak{t}_\lambda = \mu$, i.e., a tower $\mathcal{T} = \{T_\alpha : \alpha < \mu\}$. Let $\langle \alpha_\varepsilon : \varepsilon < \kappa \rangle$ be a cofinal sequence of ordinals in $\mu$. Set $\mathcal{T'} = \{T_{\alpha_\varepsilon} : \varepsilon < \kappa\}$. It is easily verified that $\mathcal{T'}$ is a tower of size $\kappa < \mu$, a contradiction.

\hfill \qedref{ttlambdareg}

We conclude this section with the following claim about the relationship between $\mathfrak{b}_\lambda$ and $\mathfrak{t}_\lambda$. It serves as an important component in the proof of the generalized Rothberger's theorem:

\begin{claim}
\label{bblamandttlam}
The relationship between $\mathfrak{b}_\lambda$ and $\mathfrak{t}_\lambda$. \newline 
Let $\lambda$ be an infinite cardinal. \newline 
\then\ $\mathfrak{t}_\lambda \leq \mathfrak{b}_\lambda$.
\end{claim}

\par \noindent \emph{Proof}. \newline 
Let $\kappa$ be any cardinal below $\mathfrak{t}_\lambda$. Assume $B = \langle b_\eta : \eta < \kappa+1 \rangle$ is a subset of ${^\lambda \lambda}$. We shall see that $B$ has an upper bound in the product ${^\lambda \lambda}$. We will choose a function $f_\eta$, by induction on $\eta \in \kappa+1$, in purpose to establish $f_\kappa$ as an upper bound.

We wish to define $f_\eta \in {^\lambda \lambda}$, strictly increasing, so that ${\rm ran}(f_\eta) \subseteq^* {\rm ran}(f_\xi)$ for every $\xi \in \eta$. This is done in the following way. Arriving at $\eta$ we choose a pseudointersection $A$ of cardinality $\lambda$ for the family $\{{\rm ran}(f_\xi) : \xi < \eta\}$. The existence of $A$ is ensured by the fact that $\eta < \mathfrak{t}_\lambda$. Now, for every $\alpha < \lambda$ we define:

$$
f_\eta(\alpha) = {\rm min} \{\alpha_* \in A : b_\eta(\beta) < \alpha_*, \forall \beta \leq \alpha\}
$$

Having the sequence $\langle f_\eta : \eta \in \kappa+1 \rangle$, we must show that $f_\kappa$ is an upper bound for $B$. For this, pick any ordinal $\xi \in \kappa$. Notice that every $f$, including $f_\kappa$, is strictly increasing (hence one to one). 

Since ${\rm ran}(f_\kappa) \subseteq^* {\rm ran}(f_\xi)$, one can choose an ordinal $\beta < \lambda$ such that $f_\kappa(\beta + \varepsilon) \in {\rm ran}(f_\xi)$ for every $\varepsilon < \lambda$. Sicne both $f_\kappa$ and $f_\xi$ are strictly increasing, we have $f_\xi(\varepsilon) \leq f_\kappa(\beta + \varepsilon)$ for every $\varepsilon < \lambda$, so for large enough $\varepsilon$ we get $f_\kappa(\varepsilon) \geq f_\xi(\varepsilon) \geq b_\xi(\varepsilon)$ as required.

\hfill \qedref{bblamandttlam}

\newpage

\section{the generalized rothberger's theorem}

Suppose $\mathcal{B}, \mathcal{C}$ are families of sets from $[\lambda]^\lambda$ (or sets from $[\lambda \times \lambda]^\lambda$). We shall say that $\mathcal{B} \perp \mathcal{C}$ if $|B \cap C| < \lambda$ for every $B \in \mathcal{B}$ and every $C \in \mathcal{C}$. A set $S$ separates $\mathcal{B}$ and $\mathcal{C}$ if $B \subseteq^* S$ for every $B \in \mathcal{B}$ and $|C \cap S| < \lambda$ for every $C \in \mathcal{C}$. This definition is symmetric (i.e., $\lambda \setminus S$ separates $\mathcal{C}$ and $\mathcal{B}$ whenever $S$ separates $\mathcal{B}$ and $\mathcal{C}$). This simple ovservation will serve us in trying to implement Rothberger's argument over singular cardinals. If there exists such $S$, we say that $\mathcal{B}$ and $\mathcal{C}$ are separable.

Following the footsteps of \cite{MR776622}, we phrase some definitions of the bounding number and prove that all are equal. We start with the fact that being unbounded on $\lambda$ and being unbounded on any member of $[\lambda]^\lambda$ are the same thing:

\begin{definition}
\label{bbb2}
The inherited bounding number, $\mathfrak{b}_\lambda^2$. \newline 
$\mathfrak{b}_\lambda^2 = {\rm min} \{|B| : B \subseteq {^\lambda \lambda}, B$ is unbounded on $X$ for every $X \in [\lambda]^\lambda\}$
\end{definition}

\begin{proposition}
\label{bbb2isbbb1}
$\mathfrak{b}_\lambda^2 = \mathfrak{b}_\lambda^1$
\end{proposition}

\par \noindent \emph{Proof}. \newline 
By the very definition of $\mathfrak{b}_\lambda^2$ it is clear that $\mathfrak{b}_\lambda^2 \geq \mathfrak{b}_\lambda$, hence $\mathfrak{b}_\lambda^2 \geq \mathfrak{b}_\lambda^1$ (due to lemma \ref{bbbisbbb1}). We shall prove that $\mathfrak{b}_\lambda^2 \leq \mathfrak{b}_\lambda^1$. Choose any $B$ which attests $\mathfrak{b}_\lambda^1$, and recall that the members of $B$ are non-decreasing functions.

We claim that for every unbounded $B \subseteq {^\lambda \lambda}$ which consists of non-decreasing functions, and every $X \in [\lambda]^\lambda$, $B$ is unbounded on $X$. Given $X$ we define for every $f \in {^\lambda \lambda}$ the upgrading function $\hat{f} \in {^\lambda \lambda}$, as follows. For every $\alpha < \lambda$ set:

$$
\hat{f}(\alpha) = {\rm sup}\{f(\gamma) : \gamma \leq {\rm min}\{\beta \in X : \beta \geq \alpha\}\}
$$

By the regularity of $\lambda$ we know that $\hat{f}$ is well defined. Now let $f$ be any function in ${^\lambda \lambda}$, and create $\hat{f}$ out of it. Choose $g \in B$ so that $\neg (g \leq^* \hat{f})$. It means that the set $J = \{\varepsilon < \lambda : \hat{f}(\varepsilon) < g(\varepsilon)\}$ is a set of size $\lambda$. For every $\varepsilon \in J$ one can pick the first ordinal $\beta_\varepsilon$ so that $\varepsilon \leq \beta_\varepsilon$ and $\beta_\varepsilon \in X$. Set $I = \{\beta_\varepsilon \in X : \varepsilon \in J\}$.

We claim that $g(\beta_\varepsilon) > f(\beta_\varepsilon)$ for every $\beta_\varepsilon \in I$, and since $I \in [X]^\lambda$ this suffices. Fix any ordinal $\varepsilon \in J$. We have $g(\beta_\varepsilon) \geq g(\varepsilon)$ (since $g$ is non-decreasing) $> \hat{f}(\varepsilon)$ (since $\varepsilon \in J$) $\geq f(\beta_\varepsilon)$ (by the definition of $\hat{f}$), so we are done.

\hfill \qedref{bbb2isbbb1}

We need some additional notation. For every $f \in {^\lambda \lambda}$ we set $L_f = \{(\alpha, \beta) \in \lambda \times \lambda : \beta \leq f(\alpha)\}$. Thinking of $f$ as a graph of a function, $L_f$ is the part 'below' $f$. Easily, $f <^* g \Leftrightarrow L_f \subset^* L_g$.

We denote the set of vertical lines in $\lambda \times \lambda$ by $\mathcal{V}_\lambda$, so $V \in \mathcal{V}_\lambda$ if there is some $\alpha < \lambda$ such that $V = \{(\alpha, \beta) : \beta < \lambda\}$. For every $X \in [\lambda \times \lambda]^\lambda$ we define $K_X = \{\alpha < \lambda : X \cap \{\alpha\} \times \lambda \neq \emptyset\}$. Now any $X \perp \mathcal{V}_\lambda$ corresponds to a function $f_X \in {^\lambda \lambda}$ which is defined as follows. $f_X(\alpha) = 0$ whenever $X \cap \{\alpha\} \times \lambda = \emptyset$ and $f_X(\alpha) = {\rm sup}\{\beta < \lambda : (\alpha,\beta) \in K_X\}$ in all other cases. Notice that this is well defined, since $\lambda$ is a regular cardinal. Notice also that $X \subseteq L_{f_X} \cap (K_X \times \lambda)$.

The following two versions of $\mathfrak{b}_\lambda$ are equal (the first is defined on $\lambda$, the second on $[\lambda \times \lambda]$, but the definitions are identical). We shall prove that the second one (hence, the first) equals $\mathfrak{b}_\lambda$:

\begin{definition}
\label{bbb5}
The bounding number via disjoint families, $\mathfrak{b}_\lambda^5$.
\begin{enumerate}
\item [$(\aleph)$] $\mathfrak{b}_\lambda^5 = {\rm min}\{|\mathcal{B}| : \mathcal{B} \subseteq [\lambda]^\lambda, \exists \mathcal{C} \subseteq [\lambda]^\lambda, |\mathcal{C}| = \lambda, \mathcal{B} \cap \mathcal{C} = \emptyset, \mathcal{B} \cup \mathcal{C}$ is almost disjoint, and for every $\mathcal{D} \in [\mathcal{C}]^\lambda$ one cannot separate $\mathcal{B}$ and $\mathcal{D}\}$
\item [$(\beth)$] $(\mathfrak{b}_\lambda^5)' = {\rm min}\{|\mathcal{B}| : \mathcal{B} \subseteq [\lambda \times \lambda]^\lambda, \exists \mathcal{C} \subseteq [\lambda \times \lambda]^\lambda, |\mathcal{C}| = \lambda, \mathcal{B} \cap \mathcal{C} = \emptyset, \mathcal{B} \cup \mathcal{C}$ is almost disjoint, and for every $\mathcal{D} \in [\mathcal{C}]^\lambda$ one cannot separate between $\mathcal{B}$ and $\mathcal{D}\}$
\end{enumerate}
\end{definition}

\begin{claim}
\label{bbb1geqbbb5}
$(\mathfrak{b}_\lambda^5)' \leq \mathfrak{b}^1_\lambda$
\end{claim}

Let $B$ be an evidence to $\mathfrak{b}^1_\lambda$, and we shall try to convert it to an evidence for $(\mathfrak{b}_\lambda^5)'$. We define $\mathcal{B} = \{L_f : f \in B\}, \mathcal{C} = \mathcal{V}_\lambda$, and almost all the requirements in the definition of $(\mathfrak{b}_\lambda^5)'$ are satisfied. It remains to prove the separating property.

Assume $\mathcal{D} \in [\mathcal{C}]^\lambda$, i.e., $\mathcal{D}$ is a collection of $\lambda$ vertical lines in $\lambda \times \lambda$. Let $D$ be the set $\{\alpha < \lambda : \{\alpha\} \times \lambda \in \mathcal{D}\}$. Let $S$ be a candidate for separating $\mathcal{B}$ from $\mathcal{D}$, i.e., $C \subseteq^* S$ for every $C \in \mathcal{D}$. We shall prove that there exits a member of $\mathcal{B}$ whose intersection with $S$ has cardinality $\lambda$.

We try to define a function $g \in {^\lambda \lambda}$ as follows. $g(\alpha) = 0$ whenever $\alpha \notin D$, and $g(\alpha) = {\rm sup}\{ \gamma < \lambda : (\alpha,\gamma) \notin S\}$ whenever $\alpha \in D$. Recall that $|D| = \lambda$, so $B$ is unbounded on $D$ (by \ref{bbb2isbbb1}). Pick $h \in B$ so that the set $I = \{\gamma \in D : h(\gamma) \geq g(\gamma)\}$ is a set of size $\lambda$. Since $h \upharpoonright I \subseteq S$ we have $|S \cap L_h| = \lambda$ hence $S$ does not separate between $\mathcal{B}$ and $\mathcal{D}$ so we are done.

\hfill \qedref{bbb1geqbbb5}

The last version of $\mathfrak{b}_\lambda$ that we need is the following:

\begin{definition}
\label{bbb7}
The non-separable bounding number, $\mathfrak{b}_\lambda^7$.
\begin{enumerate}
\item [$(\aleph)$] $\mathfrak{b}_\lambda^7 = {\rm min}\{|\mathcal{B}| : \mathcal{B} \subseteq [\lambda]^\lambda, \exists \mathcal{C} \subseteq [\lambda]^\lambda, |\mathcal{C}| = \lambda, \mathcal{B} \perp \mathcal{C}$ and one cannot separate $\mathcal{B}$ and $\mathcal{C}\}$
\item [$(\beth)$] $(\mathfrak{b}_\lambda^7)' = {\rm min}\{|\mathcal{B}| : \mathcal{B} \subseteq [\lambda \times \lambda]^\lambda, \exists \mathcal{C} \subseteq [\lambda \times \lambda]^\lambda, |\mathcal{C}| = \lambda, \mathcal{B} \perp \mathcal{C}$ and one cannot separate between $\mathcal{B}$ and $\mathcal{C}\}$
\end{enumerate}
\end{definition}

As in the case of $\mathfrak{b}_\lambda^5$, it is clear that $\mathfrak{b}_\lambda^7 = (\mathfrak{b}_\lambda^7)'$. The inequality $\mathfrak{b}_\lambda^7 \leq \mathfrak{b}_\lambda^5$ follows from the definition, so $\mathfrak{b}_\lambda = \mathfrak{b}_\lambda^1 \geq \mathfrak{b}_\lambda^5 \geq \mathfrak{b}_\lambda^7$ (by \ref{bbb1geqbbb5}). It remains to prove the last inequality that we need:

\begin{claim}
\label{bbb7isbbb}
$\mathfrak{b}_\lambda \leq \mathfrak{b}_\lambda^7$
\end{claim}

\par \noindent \emph{Proof}. \newline 
We start with $\mathcal{B}, \mathcal{C}$ as guaranteed in the definition of $(\mathfrak{b}_\lambda^7)'$. Define $\mathcal{A} = \{\bigcup \mathcal{F} : \mathcal{F} \in [\mathcal{C}]^{<\lambda} \setminus \{\emptyset\}\}$. Recall that $\mathcal{B}, \mathcal{C}$ are not separable, so in particular no member of $\mathcal{A}$ separates them. We shall construct a function $\alpha : \lambda \rightarrow \mathcal{A}$ (here we use the fact that $\lambda = \lambda^{< \lambda}$) with the following properties:

\begin{enumerate}
\item [$(a)$] $\beta < \gamma \Rightarrow \alpha(\beta) \subseteq \alpha(\gamma)$
\item [$(b)$] $\beta < \gamma \Rightarrow |\alpha(\gamma) \setminus \alpha(\beta)| = \lambda$
\item [$(c)$] For every $C \in \mathcal{C}$ there exits $\beta < \lambda$ such that $C \subseteq \alpha(\beta)$
\end{enumerate}

We do this by induction on $\lambda$. Enumerate the members of $\mathcal{C}$ by $\{C_\zeta : \zeta < \lambda \}$. For $\varepsilon = 0$ set $\alpha(0) = C_0$. For a successor ordinal $\varepsilon+1$ we have chosen some $\mathcal{F}$ such that $\alpha(\varepsilon) = \bigcup(\mathcal{F})$. Choose a member of $\mathcal{C}$, say $C'$, so that $|C' \cap \alpha(\varepsilon)| < \lambda$ (exists, by the non-separability) and define $\mathcal{F}^+ = \mathcal{F} \cup \{C'\}$. Now dedine $\alpha(\varepsilon+1) = \bigcup \mathcal{F}^+$. If $\varepsilon$ is a limit ordinal, take $\alpha(\varepsilon)$ as the union of all the previous stages. We also add to this construction $C_\varepsilon$ in the $\varepsilon$-th stage, if it is not inside yet.

By the regularity of $\lambda$ and the fact that $\mathcal{F} \in [\mathcal{C}]^{<\lambda}$, one can verify that ${\rm ran}(\alpha) \perp \mathcal{B}$ and $\mathcal{B}, {\rm ran}(\alpha)$ are not separable. So without loss of generality, $\mathcal{C} = {\rm ran}(\alpha)$. We may assume, again without loss of generality, that $\alpha(\varepsilon) = (\varepsilon+1)\times\lambda$ for every $\varepsilon < \lambda$, and then $\mathcal{B}$ and $\mathcal{V}_\lambda$ are not separable and $\mathcal{B} \perp \mathcal{V}_\lambda$.

In particular, $f_X$ is a well defined function in ${^\lambda \lambda}$ for every $X \in \mathcal{B}$. We collect these functions, stipulating $\bar{B} = \{f_X : X \in \mathcal{B}\}$, and then translating it to sets by defining $\mathcal{B'} = \{L_f : f \in \bar{B}\}$. Clearly, $\mathcal{B'} \perp \mathcal{V}_\lambda$ and $\mathcal{B'}, \mathcal{C}$ are not separable. It suffices to prove that $\bar{B}$ is unbounded in ${^\lambda \lambda}$.

Let $f \in {^\lambda \lambda}$ be any function. It follows that $L_f \perp \mathcal{V}_\lambda$. But $\mathcal{B'} \perp \mathcal{V}_\lambda$, so there must be some $g \in {^\lambda \lambda}$ such that $\neg(L_g\subseteq^*L_f)$. This means that $\neg(g\leq^*f)$, and since $g \in \bar{B}$ (because $L_g \in \mathcal{B'}$) we have shown that $\bar{B}$ is unbounded in ${^\lambda \lambda}$ and the proof is complete.

\hfill \qedref{bbb7isbbb}

Everything is ready now for the main theorem of this section:

\begin{theorem}
\label{ggrothberger}
The generalized Rothberger's theorem. \newline 
Assume $\lambda = \lambda^{<\lambda}$. \newline 
\Then\ $\mathfrak{p}_\lambda = \lambda^+ \Rightarrow \mathfrak{t}_\lambda = \lambda^+$.
\end{theorem}

\par \noindent \emph{Proof}. \newline 
If $\mathfrak{b}_\lambda = \lambda^+$ then there is nothing to prove (recall \ref{bblamandttlam}). Assume $\mathfrak{b}_\lambda > \lambda^+$. Let $\mathfrak{U} = \{u_\eta : \eta < \lambda^+\}$ be a $\mathfrak{p}_\lambda$-family. We shall construct a function $T : \lambda^+ \rightarrow [\lambda]^\lambda$ with the following two properties:

\begin{enumerate}
\item [$(a)$] $\xi < \eta < \lambda^+ \Rightarrow T_\eta \subseteq^* T_\xi$
\item [$(b)$] $u \in \mathfrak{U} \Rightarrow \exists \eta \in \lambda, T_\eta \subseteq^* u$
\end{enumerate}

If we succeed, then $\mathcal{T} = \{T_\eta : \eta < \lambda^+\}$ will be a tower (it has no pseudointersection of size $\lambda$ since such a creature will serve also as a pseudointersection for $\mathfrak{U}$). Before starting the induction process, we choose $A_\eta$ as a pseudointersection of size $\lambda$ for the collection $\{u_\xi : \xi \in \eta\}$, for every $\eta \in \lambda^+$ (possible, since the cardinality of this collection is below $\mathfrak{p}_\lambda$).

Now choose $T_\eta \in [\lambda]^\lambda$, by induction on $\lambda^+$, so that:
\begin{enumerate}
\item [$(\aleph)$] $T_\eta \subseteq^* T_\xi, u_\xi$ for every $\xi < \eta$
\item [$(\beth)$] $A_\xi \subseteq^* T_\eta$ for every $\xi \in (\eta, \lambda^+)$
\end{enumerate}

Notice that $\mathcal{T} = \{T_\eta : \eta < \lambda^+\}$ is a tower (if we can hold the induction process). In particular, it has the strong intersection property. If $\{T_\eta : \eta \in S\}$ and $|S| < \lambda$ then choose any $\gamma \in \lambda^+ \setminus {\rm sup}(S)$. $A_\gamma \subseteq^* T_\eta$ for every $\eta \in S$, hence $|A_\gamma \setminus T_\eta| < \lambda$ for every $\eta \in S$. Set $A' = \bigcup\{A_\gamma \setminus T_\eta : \eta \in S\}$ and conclude that $|A'| < \lambda$ (from the regularity of $\lambda$), hence $|A_\gamma \setminus A'| = \lambda$. Now observe that $A_\gamma \setminus A' \subseteq \bigcap\{T_\eta : \eta \in S\}$.

How do we choose the $T_\eta$-s? Set $T_0 = u_0$. Assume $\gamma > 0$ and $T_\eta$ was chosen for every $\eta < \gamma$. We define $\mathcal{B} = \{A_\xi : \xi \in (\gamma, \lambda^+)\}$ and $\mathcal{C} = \{\lambda \setminus u_\gamma\} \cup \{\lambda \setminus T_\eta : \eta \in \gamma\}$. By $(\aleph), (\beth)$ above we have $\mathcal{B} \perp \mathcal{C}$. Clearly, $|\mathcal{C}| \leq \lambda$, and by the present assumptions on $\mathfrak{b}_\lambda$ we can pick $T_\gamma \in [\lambda]^\lambda$ which separates $\mathcal{B}$ from $\mathcal{C}$.

It means that $B \subseteq^* T_\gamma$ for every $B \in \mathcal{B}$, hence requirement $(\beth)$ is satisfied. It means also that $T_\gamma \subseteq^* \lambda \setminus C$ for every $C \in \mathcal{C}$, which gives requirement $(\aleph)$, so we are done.

\hfill \qedref{ggrothberger}

\begin{corollary}
\label{cofcontinuum}
Assume $\lambda = \lambda^{< \lambda}$. \newline 
If $\cf(2^\lambda) \in \{\lambda^+,\lambda^{++}\}$, \then\ $\mathfrak{p}_\lambda = \mathfrak{t}_\lambda$.
\end{corollary}

\par \noindent \emph{Proof}. \newline 
As in the case of $\lambda = \omega$, we have $\lambda \leq \kappa < \mathfrak{t}_\lambda \Rightarrow 2^\kappa = 2^\lambda$. Since $2^{\cf(2^\lambda)} > \lambda$ we infer that $\mathfrak{t}_\lambda \leq \cf(2^\lambda)$, and the corollary follows.

\hfill \qedref{cofcontinuum}

\newpage 

\section{conclusions}

We shall derive, in this section, some conclusions from the results of the former section. The first one is about the cofinality of $\mathfrak{p}_\lambda$, and the second is about the relationship between $\mathfrak{p}_\lambda$ and $\mathfrak{t}_\lambda$. The idea is that the main argument in the proof of Rothberger's theorem applies to singulars with low cofinality. So if $\mathfrak{p}_\lambda$ is such a cardinal, then $\mathfrak{t}_\lambda$ is also a singular cardinal, which is impossible.

\begin{claim}
\label{ccofpilam}
The cofinality of $\mathfrak{p}_\lambda$. \newline 
Assume $\lambda = \lambda^{<\lambda}$. \Then\ $\cf(\mathfrak{p}_\lambda) \neq \lambda$.
\end{claim}

\par \noindent \emph{Proof}. \newline 
Before we start, notice that $\mathfrak{p}_\lambda > \lambda$. Indeed, if $\mathcal{F}$ has the strong intersection property and $|\mathcal{F}| \leq \lambda$, then it is easy to construct a pseudo intersection of size $\lambda$ for $\mathcal{F}$ by diagonalizing.

Assume $\cf(\mu) = \lambda < \mu$. Suppose toward contradiction that $\mathfrak{p}_\lambda = \mu$. Since $\mathfrak{p}_\lambda \leq \mathfrak{t}_\lambda \leq \mathfrak{b}_\lambda$ and $\mathfrak{t}_\lambda$ is a regular cardinal, we know that $\mathfrak{t}_\lambda > \mu$ hence $\mathfrak{b}_\lambda > \mu$.
Let $\mathfrak{U} = \{u_\eta : \eta < \mu\}$ be an evidence for $\mathfrak{p}_\lambda$. We shall produce a tower $\mathfrak{T}$ out of $\mathfrak{U}$, keeping its cardinality.

For this, we try to define a function $T : \mu \rightarrow [\lambda]^\lambda$ so that $T_\eta \subseteq^* T_\xi$ for every $\xi < \eta < \mu$ and for every $u \in \mathfrak{U}$ there exists $\eta < \mu$ so that $T_\eta \subseteq^* u$. If we succeed, then $\mathfrak{T} = \{T_\eta : \eta < \mu\}$ will be a tower. The fact that for every $u \in \mathfrak{U}$ there exists $\eta < \mu$ so that $T_\eta \subseteq^* u$ with the fact that $\mathfrak{U}$ is an evidence for $\mathfrak{p}_\lambda$ ensures that $\mathfrak{T}$ has no pseudointersection of size $\lambda$.

At first stage, Choose a cofinal sequence $\langle \alpha_\varepsilon : \varepsilon < \cf(\mu) \rangle$.
Let $A_{\alpha_\varepsilon}$ be a pseudointersection of size $\lambda$ of $\{u_\xi : \xi < \alpha_\varepsilon\}$, for every $\alpha_\varepsilon < \cf(\mu)$. Now we choose, by induction on $\eta < \mu$, a set $T_\eta \in [\lambda]^\lambda$ such that two things are maintained along the way:

\begin{enumerate}
\item [$(a)$] $T_\eta \subseteq^* T_\xi, u_\xi$ for every $\xi < \eta$
\item [$(b)$] $A_{\alpha_\varepsilon} \subseteq^* T_\eta$ for every $\alpha_\varepsilon > \eta$
\end{enumerate}

We start, in the stage of $\gamma = 0$, with $T_0 = u_0$. Arriving at $\gamma > 0$ and assuming that $T_\eta$ was chosen for every $\eta < \gamma$, we define $\mathcal{B} = \{A_{\alpha_\varepsilon} : \gamma< \alpha_\varepsilon, \varepsilon < \cf(\mu)\}$ and $\mathcal{C} = \{\lambda \setminus u_\gamma\} \cup \{ \lambda \setminus T_\eta : \eta < \gamma\}$. Notice that $|\mathcal{B}| \leq \cf(\mu) = \lambda$ and $|\mathcal{C}| < \mu < \mathfrak{b}_\lambda$.

We claim that $\mathcal{B} \perp \mathcal{C}$. Indeed, pick any $A_{\alpha_\varepsilon} \in \mathcal{B}$. Since $\gamma < \alpha_\varepsilon$ we have $A_{\alpha_\varepsilon} \subseteq^* u_\gamma$ (by the definition of $A_{\alpha_\varepsilon}$) hence $|A_{\alpha_\varepsilon} \cap \lambda \setminus u_\gamma| < \lambda$. Similarly, for every $\eta < \gamma$ we have $A_{\alpha_\varepsilon} \subseteq^* T_\eta$ (by the induction hypothesis) so $|A_{\alpha_\varepsilon} \cap \lambda \setminus T_\eta| < \lambda$ as before.

We can use now the fact that $|\mathcal{C}| < \mathfrak{b}_\lambda^7$ and $|\mathcal{B}| \leq \lambda$ to choose a set $T_\gamma \in [\lambda]^\lambda$ which separates $\mathcal{B}$ from $\mathcal{C}$. So $A_{\alpha_\varepsilon} \subseteq^* T_\gamma$ for every $\gamma < \alpha_\varepsilon$ and this is requirement $(b)$ in the induction process. On the other hand, $T_\gamma \subseteq^* \lambda \setminus C$ for every $C \in \mathcal{C}$ which means that $T_\gamma \subseteq^* T_\xi$ for every $\xi < \gamma$, and $T_\gamma \subseteq^* u_\gamma$, and that gives also $T_\gamma \subseteq^* T_\xi \subseteq^* u_\xi$ for every $\xi < \gamma$. These relations accomplish the construction.

We claim now that $\mathfrak{T} = \{T_\eta : \eta < \mu\}$ is a tower. Property $(a)$ ensures that it has no pseudo intersection of size $\lambda$, since such an example will be also a pseudo intersection of the $\mathfrak{p}_\lambda$ family. Property $(b)$ makes sure that $\mathfrak{T}$ has the strong intersection property. But there is no tower whose cardinality is $\mu$, so the proof is complete.

\hfill \qedref{ccofpilam}

Our next mission is to say something about the distance between $\mathfrak{p}_\lambda$ and $\mathfrak{t}_\lambda$. The main point here is that Rothberger's theorem converts any $\mathfrak{p}_\lambda$-family into a $\mathfrak{t}_\lambda$-family of the same size, even when the cardinality of these families is above $\mathfrak{p}_\lambda$. Let us start with the following definition:

\begin{definition}
\label{ppfamily}
Irredundant $\mathfrak{p}_\lambda$-families. \newline 
Suppose $\lambda = \lambda^{<\lambda}$, and $\theta \leq \lambda$.
\begin{enumerate}
\item [$(\alpha)$] The collection $\mathcal{U} = \{u_\eta : \eta < \kappa\}$ is a $\mathfrak{p}_\lambda$-family if $\mathcal{U}$ has the strong intersection property and there is no pseudointersection of size $\lambda$ for $\mathcal{U}$ (notice that $\kappa$ might be greater than $\lambda$).
\item [$(\beta)$] A $\mathfrak{p}_\lambda$-family $\mathcal{U}$ is $\theta$-irredundant if there exits $\mathcal{U'} \subseteq \mathcal{U}$, $|\mathcal{U'}| = \theta$, such that $\mathcal{U} \setminus \mathcal{U'}$ is not a $\mathfrak{p}_\lambda$-family. Moreover, removing from $\mathcal{U'}$ a bounded subset does not change this fact.
\item [$(\gamma)$] A $\mathfrak{p}_\lambda$-family $\mathcal{U} = \{u_\eta : \eta < \kappa\}$ is accurate, if $\{u_\eta : \eta < \delta\}$ has pseudointersection of size $\lambda$ for every $\delta < \kappa$.
\end{enumerate}
\end{definition}

\begin{lemma}
\label{bbigpfamily}
Assume $\mu > \cf(\mu) = \theta$. \newline 
Suppose $\mathfrak{p}_\lambda < \mu < \mathfrak{t}_\lambda$, and $\mathcal{U}$ is a $\theta$-irredundant $\mathfrak{p}_\lambda$-family. \newline 
\then\ there exits an accurate $\mathfrak{p}_\lambda$-family of size $\mu$.
\end{lemma}

\par \noindent \emph{Proof}. \newline 
For $\mathfrak{p}_\lambda$ itself there is an accurate $\mathfrak{p}_\lambda$-family, so let $\mathcal{U} = \{u_\eta : \eta < \mathfrak{p}_\lambda\}$ be an accurate $\mathfrak{p}_\lambda$-family, and let $\mathcal{T} = \{T_\beta : \beta < \mathfrak{t}_\lambda\}$ exemplify $\mathfrak{t}_\lambda$. We cut the tail of the tower, so $\mathcal{T'} = \{T_\beta : \beta < \mu\}$ has a pseudointersection of size $\lambda$, call it $A$ (recall that $\mu < \mathfrak{t}_\lambda$). Choose $\mathcal{U'} \subseteq \mathcal{U}, |\mathcal{U'}| = \theta$ so that $\mathcal{U} \setminus \mathcal{U'}$ is not a $\mathfrak{p}_\lambda$-family.

We may assume that every member of $\mathcal{U}$ is contained in $A$ (use a bijection $h : \lambda \rightarrow A$, to transmit the members of $\mathcal{U}$ into subsets of $A$). Now set $\mathcal{R} = \mathcal{T'} \cup \mathcal{U}$. This is a $\mathfrak{p}_\lambda$-family of size $\mu$. Enumerate its members in such a way that the members of $\mathcal{U'}$ are cofinal, so this $\mathfrak{p}_\lambda$-family is accurate.

\hfill \qedref{bbigpfamily}

On the base of this lemma, we can try to say something about the difference between $\mathfrak{p}_\lambda$ and $\mathfrak{t}_\lambda$. The theorem below is applicable also for the case $\lambda = \omega$, i.e., the common cardinal invariants $\mathfrak{p}$ and $\mathfrak{t}$.

\begin{theorem}
\label{pppttt}
The discrepancy between $\mathfrak{p}_\lambda$ and $\mathfrak{t}_\lambda$. \newline 
Assume $\lambda = \lambda^{<\lambda}$. \newline 
Suppose there exits a $\lambda$-irredundant $\mathcal{U}$ which exemplifies $\mathfrak{p}_\lambda$. \newline 
\Then\ there is no singular cardinal $\mu$ of cofinality $\lambda$ between $\mathfrak{p}_\lambda$ and $\mathfrak{t}_\lambda$. \newline 
In particular, if there exits an $\aleph_0$-irredundant $\mathfrak{p}$-family, then $\mathfrak{t} \leq \mathfrak{p}^{+n}$ for some $n \in \omega$.
\end{theorem}

\par \noindent \emph{Proof}. \newline 
Assume toward contradiction that $\mathfrak{p}_\lambda < \mu < \mathfrak{t}_\lambda$, and $\mu > \cf(\mu) = \lambda$. By lemma \ref{bbigpfamily} we can choose an accurate $\mathfrak{p}_\lambda$-family $\mathfrak{U} = \{u_\eta : \eta < \mu\}$. Let $\mathfrak{U'}$ be an evidence to the fact that $\mathfrak{U}$ is $\lambda$-irredundant. We may assume that the members of $\mathfrak{U'}$ are cofinal in $\mathfrak{U}$. As in the beginning of this section, we shall produce a tower of size $\mu$ out of $\mathfrak{U}$, which is absurd.

We choose a cofinal sequence in $\mu$, namely $\langle \alpha_\varepsilon : \varepsilon < \cf(\mu) \rangle$. We fix a sequence of sets $\langle A_{\alpha_\varepsilon} : \varepsilon < \cf(\mu) \rangle$ whence $A_{\alpha_\varepsilon}$ is a pseudointersection of size $\lambda$ for $\{u_\xi : \xi < \alpha_\varepsilon\}$ (the existence of $A_{\alpha_\varepsilon}$ is justified by the fact that $\mathfrak{U}$ is accurate).

Now we define $T_\eta \in [\lambda]^\lambda$ by induction on $\eta < \mu$ so that:

\begin{enumerate}
\item [$(a)$] $T_\eta \subseteq^* T_\xi, u_\xi$ for every $\xi < \eta$
\item [$(b)$] $A_{\alpha_\varepsilon} \subseteq^* T_\eta$ for every $\alpha_\varepsilon > \eta$
\end{enumerate}

As before, $T_0 = u_0$. Assume $\gamma > 0$ and $T_\eta$ was chosen for every $\eta < \gamma$. Aiming to choose $T_\gamma$, set $\mathcal{B} = \{A_{\alpha_\varepsilon} : \gamma< \alpha_\varepsilon, \varepsilon < \cf(\mu)\}$ and $\mathcal{C} = \{\lambda \setminus u_\gamma\} \cup \{ \lambda \setminus T_\eta : \eta < \gamma\}$. Clearly, $|\mathcal{B}| \leq \lambda$ and $|\mathcal{C}| < \mathfrak{b}_\lambda$.

As above, $\mathcal{B} \perp \mathcal{C}$, so choose $T_\gamma \in [\lambda]^\lambda$ which separates $\mathcal{B}$ and $\mathcal{C}$. It means that $A_{\alpha_\varepsilon} \subseteq^* T_\gamma$ whenever $\gamma < \alpha_\varepsilon$, so requirement $(b)$ is satisfied. It also means that $T_\gamma \subseteq^* \lambda \setminus C$ for every $C \in \mathcal{C}$ which gives the requirements in $(a)$. Having the induction hypothesis, we can create the collection $\mathcal{T} = \{T_\eta : \eta < \mu\}$.

We claim that $\mathcal{T}$ is a tower. The verification of all the properties of a tower is as in the previous claims. In particular, $\mathcal{T}$ has the strong intersection property. If $S \subseteq \mu, |S| < \lambda$ then $|S|$ is bounded in $\mu$. Choose any $\alpha_\varepsilon$ above ${\rm sup}(S)$. It follows that $A_{\alpha_\varepsilon} \subseteq^* T_\eta$ for every $\eta \in S$. Define $A' = \bigcup\{A_{\alpha_\varepsilon} \setminus T_\eta : \eta \in S\}$. Clearly, $|A'| < \lambda$. Now $A_{\alpha_\varepsilon} \setminus A' \subseteq \bigcap\{T_\eta : \eta \in S\}$.

As before, the existence of a tower of size $\mu$ is an absurd, so the proof is complete.

\hfill \qedref{pppttt}

It is worth to rephrase the theorem above for the specific case of $\lambda = \omega$ in the following way:

\begin{corollary}
\label{ppeqttt}
Suppose there exits an $\aleph_0$-irredundant $\mathfrak{p}$-family, and $\lambda > \cf(\lambda) = \aleph_0$.
\begin{enumerate}
\item [$(a)$] If $\mathfrak{t} = \lambda^+$ (e.g., $\mathfrak{t} = \aleph_{\omega+1}$) then $\mathfrak{p}=\mathfrak{t}$
\item [$(b)$] If $\mathfrak{t}$ is weakly inaccessible then $\mathfrak{p}=\mathfrak{t}$
\end{enumerate}
\end{corollary}

\newpage

\bibliographystyle{amsplain}
\bibliography{arlist}

\providecommand{\bysame}{\leavevmode\hbox to3em{\hrulefill}\thinspace}
\providecommand{\MR}{\relax\ifhmode\unskip\space\fi MR }
\providecommand{\MRhref}[2]{%
  \href{http://www.ams.org/mathscinet-getitem?mr=#1}{#2}
}
\providecommand{\href}[2]{#2}
\begin{thebibliography}{1}

\bibitem{MR0004281}
Fritz Rothberger, \emph{Sur les familles ind\'enombrables de suites de nombres
  naturels et les probl\`emes concernant la propri\'et\'e {$C$}}, Proc.
  Cambridge Philos. Soc. \textbf{37} (1941), 109--126. \MR{MR0004281 (2,352a)}

\bibitem{MR2518968}
Saharon Shelah, \emph{A comment on ``{$\mathfrak p<\mathfrak t$}''}, Canad.
  Math. Bull. \textbf{52} (2009), no.~2, 303--314. \MR{MR2518968 (2010h:03078)}

\bibitem{MR776622}
Eric~K. van Douwen, \emph{The integers and topology}, Handbook of set-theoretic
  topology, North-Holland, Amsterdam, 1984, pp.~111--167. \MR{MR776622
  (87f:54008)}

\end{thebibliography}

\end{document}